\documentclass[a4paper, 
11pt
]{article}

\usepackage{graphicx}
\usepackage{amsmath}
\usepackage{amsfonts}
\usepackage{amssymb}
\usepackage{url}
\usepackage[caption=false]{subfig}
\usepackage{multirow}
\usepackage{authblk}

\begin{document}

\title{Advances in Reduced Order Methods for Parametric Industrial Problems in Computational Fluid Dynamics}

\author[]{Gianluigi~Rozza\footnote{gianluigi.rozza@sissa.it}}
\author[]{Haris~Malik\footnote{mmalik@sissa.it}}
\author[]{Nicola~Demo\footnote{nicola.demo@sissa.it}}
\author[]{Marco~Tezzele\footnote{marco.tezzele@sissa.it}}
\author[]{Michele~Girfoglio\footnote{mgirfogl@sissa.it}}
\author[]{Giovanni~Stabile\footnote{gstabile@sissa.it}}
\author[]{Andrea~Mola\footnote{andrea.mola@sissa.it}}

\affil[]{Mathematics Area, mathLab, SISSA, International School of Advanced Studies, via Bonomea 265, I-34136 Trieste, Italy}

\maketitle

\begin{abstract}
Reduced order modeling has gained considerable attention in recent
decades owing to the advantages offered in reduced computational times
and multiple solutions for parametric problems. The focus of this
manuscript is the application of model order reduction techniques in
various engineering and scientific applications including but not
limited to mechanical, naval and aeronautical engineering. The focus
here is kept limited to computational fluid mechanics and related
applications. The advances in the reduced order modeling with proper
orthogonal decomposition and reduced basis method are presented as
well as a brief discussion of dynamic mode decomposition and also some
present advances in the parameter space reduction. Here, an overview
of the challenges faced and possible solutions are presented with
examples from various problems. 
\end{abstract}

\section{Introduction and Motivation}
\label{sec:Intro}

Advances in computational capabilities and the computational power of
modern computer systems lead to more accurate and complex mathematical
and numerical models. These ever increasing complexity is coordinated
with the progress in the modeling, numerical analysis and faster
algorithms in computational science. Despite all the advancements,
there are still several problems in science and engineering which are
relatively difficult to be computed \cite{Chinesta_Book}. The various
challenges faced during can be summarized as: 

\begin{itemize}
\item[-] The first challenge is handling of the higher dimensional
  space. Models with definition in high dimensional spaces encounter
  what is called ``curse of dimensionality", examples of such problems
  are many parametric problems.

\item[-] On-line control of complex systems, which requires fast real
  time simulations usually on hand held devices such as tablets.

\item[-] Problems involving inverse identification, process and shape
  optimization. 

\end{itemize}

This list provides very basic examples and in practical there are many
more challenges that remain out of the reach of traditional
computational strategies. There are two possible solutions to these
challenges; the first is to reduce the complexity of the mathematical
model by introducing assumptions in the physics of the problem to
simplify the mathematical model, but this approach is not always
practical and can also be the cause of introduction of errors. The
second approach is the utilization of High Performance Computing (HPC)
\cite{Malik_thesis}, which has its own drawback as HPC is expensive to
install and immobile. A more suitable approach is the development of
reduced order models, which relies on the division of the problem such
that a complex and higher order problem is solved with greater
accuracy using expensive computational facilities normally during
off-line phase and then using fast on-line phase to compute specific
solutions on much less computational expense. 

Reduced order modeling relies on the mathematical approach rather than
introducing assumptions to simplify the problem \cite{Malik_thesis},
to obtain much smaller model than the high fidelity model without
compromising the accuracy of the solution. There are a number of
techniques available for the purpose of model order reduction;
categorized into two types; \textit{a posteriori} and \textit{a
  priori}. A posteriori model order reduction method includes proper
orthogonal decomposition (POD) \cite{Berkooz, LIANG2002527},
trajectory piecewise linear method (TPWL) \cite{Rewienski, Malik2016},
reduced basis (RB) method \cite{quarteroniRB2016,
  quarteroni2007numerical, Drohmann2012} including also tools like
empirical interpolation method (EIM) \cite{BARRAULT2004667} to name a
few. A priori model order reduction methods include methods such as
proper generalized decomposition (PGD) \cite{Chinesta2011, Malik2018}
and a priori reduction method (APR) \cite{Allery2011, Ryckelynck2006}.  

In this manuscript, the authors present an overview of various
applications of reduced order modeling techniques with the focus on
computational fluid dynamics (CFD). The paper is organized as
Section~\ref{sec:Advances} provides the recent advances in the reduced
order modeling, from Section~\ref{sec:CFD_FV} to
Section~\ref{sec:Industrial} we present, application of reduced order
modeling (ROM) in several scientific and industrial applications.

\section{Advances in Reduced Order Modeling}
\label{sec:Advances}

With the recent advances in the computational science, the focus of
research is more towards the development of numerical methods and
strategies for the parametric problems involving partial differential
equations (PDEs). The introduction of parameters as discussed in
Section~\ref{sec:Intro} increases the dimension of the problem space,
the parameters can arise from material, geometrical and
non-dimensional coefficients. The research in the fields of numerical
analysis specifically in computational mechanics with applications
such as simulation, optimization and real time control deals with such
parametrized PDEs. Such cases require multiple numerical solutions of
PDEs with different parameter values which require high computational
efficiency. These problems, therefore, provide the need for
development of reduced order modeling techniques such as RB, POD and
PGD.  

Reduced basis methods \cite{hesthaven2015certified} have been
developed into a strong model order reduction method in previous years
providing reduction of computational times for the solution of
parametrized PDEs. Similar to most of model reduction techniques,
reduced basis methods divide the solution in an off-line stage and an
on-line stage. In the off-line stage, a solution is sought of the high
fidelity model; of the order, say N; with the help of suitable
discretization technique finite element method (FEM), finite volume
method (FVM) and finite difference (FD) etc., depending upon the
nature of the problem. During this stage, a number of solutions are
stored for different parameter values which are chosen in an optimal
way and are subsequently used to generate a reduced basis of much
smaller order M $\ll$ N. Once the off-line phase is completed, reduced
basis functions can be utilized for the generation of new solutions
for new parameter values combining the previously computed basis
functions by means of a Galerkin projection
\cite{hesthaven2015certified, rozza_2011} in a fast and efficient
manner. This problem in reduced dimensional space is very small and
therefore useful in the deployment of real time scenarios. The
solution thus obtained is reliable and accurate ensured by
residual-based a posteriori error estimates.  

Recent research activities in the field of reduced order modeling
prove the effectiveness of this approach and has resulted in
significant development of model order reduction methods for several
different problems of practical interest \cite{Quarteroni2009}. For
engineering applications, in addition to be able to perform efficient
numerical simulations for complex geometries of various different
materials, reduced order methods need to be capable of parameterizing
the geometric shape of the structure itself. Similar, requirements are
also found in medical applications for example the CFD analysis of
blood flow through vessels \cite{TYFA2018228, BallarinJCP}. Current
focus of the research is towards developing the theory and the
methodology of reduced order methods for computational fluid dynamics
for different physical and temporal scales and also for complex
nonlinear problems such as bifurcations and instabilities. On a vastly
different scale from the blood flows
\cite{Ballarin_thesis,BallarinJCP} and biology of singular cell and
micro-organisms \cite{Alouges2013}, is the application in the naval
and nautical engineering. As well as there are applications in well
defined CFD problems of aerospace, mechanical and automotive
engineering and porous media and geophysics.

These developed reduced order methods can be used in combination with
techniques in data assimilation and uncertainty quantification for the
solution of complex inverse problems found in the multidisciplinary
fields described here.

\section{CFD with Finite Volume Method}
\label{sec:CFD_FV}

\begin{figure}[t]
\centering
\begin{tabular}{cc}
\subfloat[Naval Engineering\label{fig:naval_engineering}]{\includegraphics[width=0.5\textwidth]{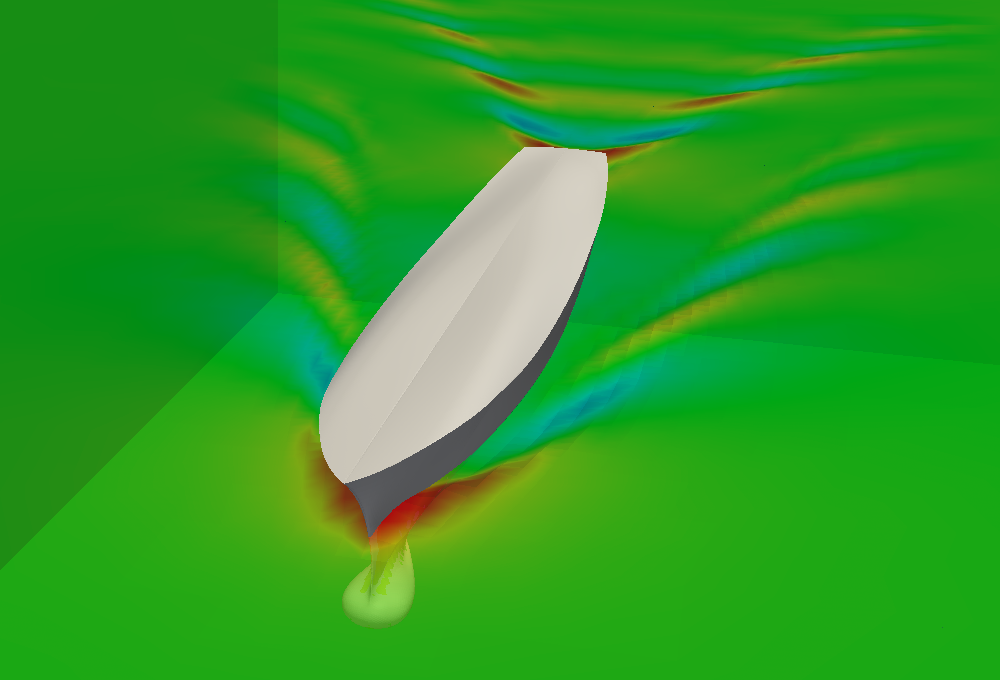}} & 
\subfloat[Aeronautics\label{fig:aeronautics}]{\includegraphics[width=0.49\textwidth]{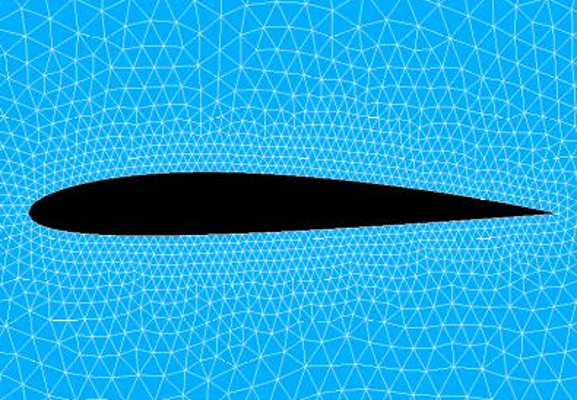}} \\
\multicolumn{2}{c}{\subfloat[Industrial Application\label{fig:industrial_application}]{\includegraphics[width=0.5\textwidth]{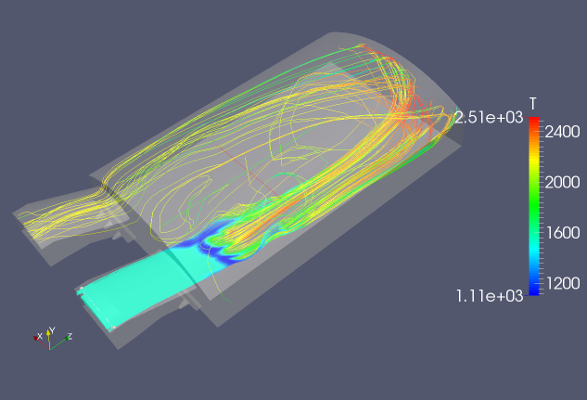}}} \\
\end{tabular}
\caption{Possible applications where the finite volume discretization is particularly widespread}
\label{fig:cfd_fv}
\end{figure}

The finite volume method is particularly widespread in several
engineering fields (aeronautics engineering, naval engineering,
automotive engineering, civil engineering, ...) and historically is
widely used in industrial applications characterized by higher values
of the Reynolds number, refer Figure~\ref{fig:cfd_fv} for typical
examples. The finite volume method rather than operating on the strong
form of the equations works on the integral version of the
equations. Divergence terms are then converted into to surface
integrals exploiting the divergence theorem and the conservation law
is enforced in each finite volume. For more details concerning the
derivation of the discretized equations and the mathematical
foundations of the method the reader may see
\cite{ferziger99:CMFD,ohlbergerFV}. However, in the field of reduced
order methods, this discretization technique is usually less employed
respect to the finite element method. To this purpose, one of the
first contribution dealing with finite volumes and the reduced basis
method, in the context of a linear evolution equations, can be found
in \cite{Haasdonk2008}. A recent contribution considering non-linear
evolution equations modeling fluid dynamics problems can be found
in~\cite{Stabile2017}, in~\cite{stabile_stabilized} pressure
stabilization techniques normally employed in a FEM-based POD-Galerkin
methods are extended to a finite volume framework
while~\cite{Carlberg2018} proposes a new ROM where conservation laws
are enforced also at the reduced order level. The implementation of
ROMs for finite volume schemes will allow to more effectively propose
the reduced order methodology outside the academic environment for
complex and real-world problems that industrial partners face on a
daily basis. Recently the open-source library ITHACA-FV~\cite{ITHACA}
has been released which is based on OpenFOAM~\cite{OF}, the most
widely used general purpose open-source CFD software package. 

\subsection{Some ROM challenges in CFD for finite volume scheme}
\label{sec:ROM_challenges}
In the last years some progresses have been achieved but there are
still several issues that need to be tackled. Among them, the most
challenging ones are ROMs for turbulent flows and ROMs that include
geometrical parametrization. 
Regardless from the starting full order discretization technique, the
majority of projection based ROMs are limited to laminar flows and
relatively few reports for turbulent flows have appeared in the
literature and here we report some of the most relevant
\cite{taddei2017, wang_turb,Lorenzi2016,Iollo2000_b,noack2005}. Among
these works only few are based on finite volumes, on the contrary, at
full order level, a large variety of closure models for turbulent
flows can be found. For this reason it is crucial to export what has
been done for full order finite volume schemes to a reduced setting. 
Concerning geometrical parametrization, in finite elements, a common
strategy  is the usage of equation written into a reference domain in
order to have all the results in a common domain. This approach is
however not easy to be transferred to a finite volume setting and in
many cases even not possible. In fact, dealing with non-linear and
non-explicit schemes, or with complex geometrical deformations, this
operation becomes not possible. Possible ways to overcome this
limitation could rely on the usage of immersed methods
\cite{stabile2018ibm}. In this way is possible to write all the
equations on the same physical domain and to treat the immersed
structure as an external forcing term. 

\section{Dynamic Mode Decomposition}
\label{sec:DMD}
\begin{figure}[t]
\centering
\includegraphics[width=11cm]{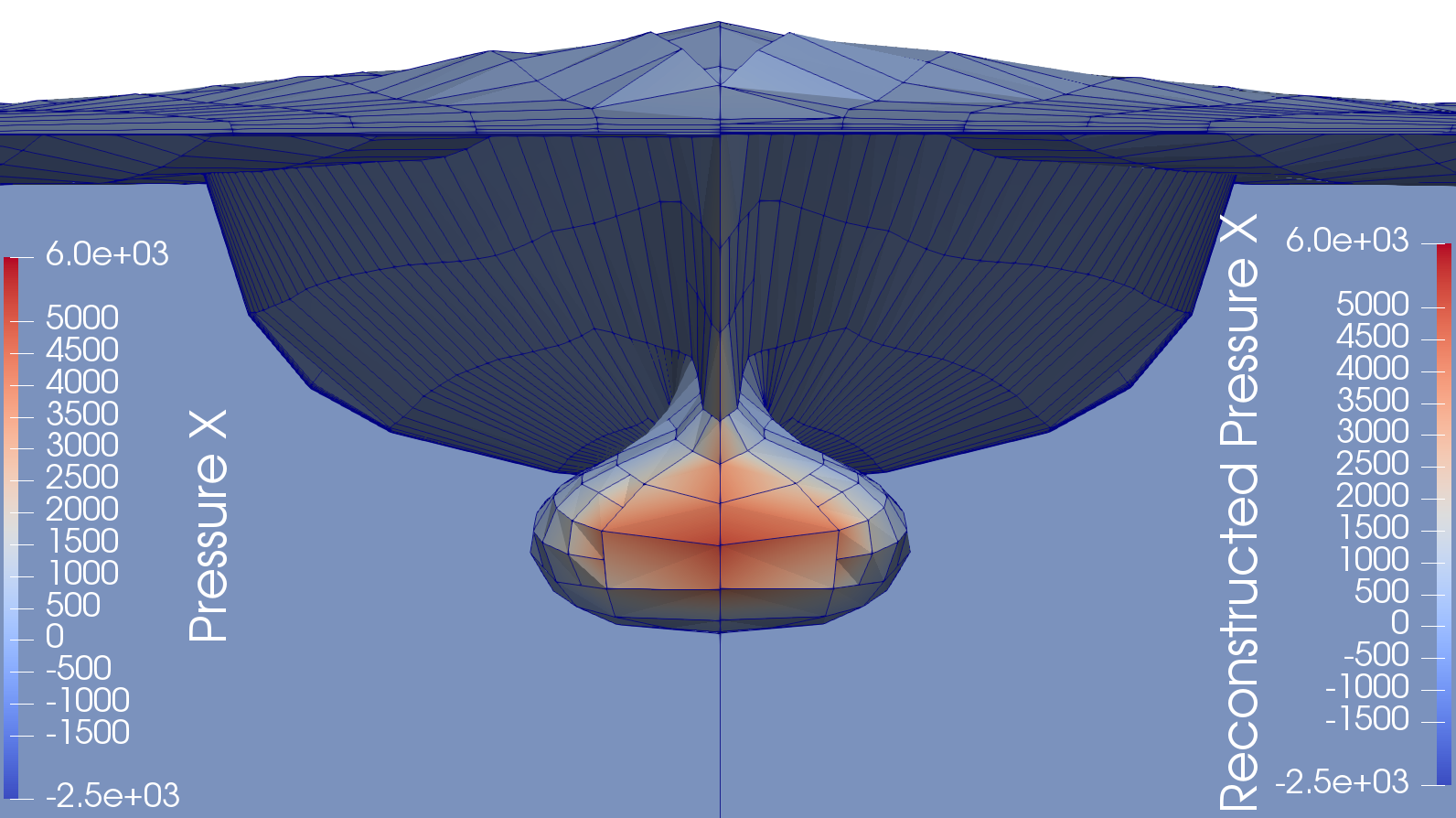}
\caption{Application of DMD on a naval simulation: on the left the
  high-fidelity solution, on the right the points displacement and the
  pressure field reconstructed.} 
\label{fig:dmd}
\end{figure}

Dynamic mode decomposition (DMD) is an emerging tool for complex
systems analysis. Initially introduced in the fluid dynamics community
\cite{schmid2010dynamic}, this technique has be adopted also in many
other fields due to the capability to represent a complex --- also
nonlinear --- system as linear combination of few main structures that
evolve linearly in time. Since DMD approximates the Koopman operator
\cite{koopman1931hamiltonian} using just the data extracted from the
underlying system, this method does not require any information about
the governing equations. 
Due to its diffusion, in the last years several variants of the
standard algorithm have been developed, like multiresolution DMD
\cite{kutz2016multiresolution}, forward-backward DMD
\cite{dawson2016characterizing}, higher order DMD \cite{le2017higher}
and compressed DMD \cite{erichson2016compressed}. All these variants
have been implemented in an open source Python package called PyDMD
\cite{demo18pydmd}. An example of the PyDMD application on a fluid
dynamics simulation is shown in Fig~\ref{fig:dmd}. 

Basically, we want to approximate the Koopman operator $\mathbf{A}$
such that it can simulates the system time evolution, hence the
relation between two sequential instants is
$\mathbf{x}_{k+1}~=~\mathbf{A}\mathbf{x}_k$, where $\mathbf{x}_i$
denotes the state at $i$-th instant. To achieve this, we collect a
series of data vectors containing the system states, from now on
\textit{snapshots}, and arrange them into two matrices such that
$\mathbf{S} = \begin{bmatrix}\mathbf{x}_1 & \dots &
  \mathbf{x}_{m-1}\end{bmatrix}$ and $\mathbf{\dot{S}}
= \begin{bmatrix}\mathbf{x}_2 & \dots &
  \mathbf{x}_{m}\end{bmatrix}$. The linear operator is built to
minimize the error $\|\mathbf{\dot{S}} - \mathbf{A}\mathbf{S}\|$, so
$\mathbf{A} = \mathbf{\dot{S}}\mathbf{S}^\dagger$, where $^\dagger$
refers to the Moore-Penrose pseudo-inverse. 
To avoid handling this large matrix, the reduced operator is
computed. The matrix $\mathbf{A}$ is projected on the spaced spanned
by the left-singular vectors of matrix $\mathbf{S}$, found by the
truncated singular value decomposition (SVD), that is $\mathbf{S}
\approx \mathbf{U}_r \mathbf{\Sigma}_r\mathbf{V}_r^*$. Once we
obtained the reduced operator $\mathbf{\tilde{A}}$, we can reconstruct
the eigenvectors and eigenvalues of the matrix $\mathbf{A}$ thanks to
the eigen-decomposition of $\mathbf{\tilde{A}}$, which is
$\mathbf{\tilde{A}} \mathbf{\Lambda} =  \mathbf{\Lambda}
\mathbf{W}$. In particular the elements in $\mathbf{\Lambda}$
correspond to the nonzero eigenvalues of $\mathbf{A}$, while the real
eigenvectors, the so called exact modes \cite{tu2014dynamic},
can be computed as $\mathbf{\Phi} = \mathbf{\dot{S}}\mathbf{V}_r
\mathbf{\Sigma}_r^{-1}\mathbf{W}$. 

\section{Proper Orthogonal Decomposition with Interpolation}
\label{sec:PODI}
Within the model reduction techniques, the proper orthogonal
decomposition interpolation (PODI) method is exploited both in the
academia and industrial context for the real time approximation of the
solution of parametric partial differential equations
\cite{bui2003proper, schilders2008model,
  wojciak2011investigation,salmoiraghi2018free, demo2018shape}. Since
it relies only on the high-fidelity solutions, its biggest advantage
is to be completely independent from the full-order solver and to not
require any assumptions on the underlying system. 

Initially, in the PODI approach, the parametric space $\mathbb{D}^N$
is sampled and the high-fidelity solutions are computed by solving the
full-order model built using these parametric points. At this stage,
the numerical method and the discretization can be chose to reach the
desired accuracy. Once these solutions are collected --- the most
expensive phase --- we can represent them as a linear combination of
few main structures, the so called POD \textit{modes}. The modes
correspond to the left-singular vectors individuated by applying the
singular value decomposition (SVD) to the solutions matrix. We define
the \textit{modal coefficients} as the projection of the high-fidelity
solutions onto the space spanned by the modes: thanks to the
correlation between the parametric points and the modal coefficients,
we can interpolate the solution for any parameter point that belongs
to $\mathbb{D}^N$. Even though the accuracy of the approximated
solution depends from the interpolation method used, the capability to
generate the reduced model using only the system output make this
equation-free algorithm specially suited for the industrial
applications. About this, the PODI method has been proposed and
implemented by the mathLab group in an open-source package on Github,
called EZyRB~\cite{ezyrb}.

\section{Efficient Geometrical Parameterization Techniques in the Context of ROM}
\label{sec:Efficient}

Nowadays shape optimization has gained a lot of interest. In this
framework an efficient and accurate geometrical parametrization is a
critical part for each optimal shape design simulation campaign. This
is true especially in the context of reduced order modeling, where it
is possible to discern the shape morphing methods in two main groups:
\emph{general purpose} or \emph{problem specific}. While the later
aims at reducing the parameter space using specific characteristics of
the problem at hand, like, for instance, the centerlines-based
approach proposed in~\cite{BallarinJCP}, the first approach applies to
a wide range of problems. 
Among possible \emph{general purpose} methods we mention \emph{Free-Form Deformation} (FFD)~\cite{sederbergparry1986,LassilaRozza2010,rozza2013free,salmoiraghi2018free},
\emph{Radial Basis Functions} (RBF)
interpolation~\cite{buhmann2003radial,morris2008cfd,manzoni2012model}
or \emph{Inverse Distance Weighting} (IDW)
interpolation~\cite{shepard1968,witteveenbijl2009,forti2014efficient}. These
methods identify the parameters as the displacements of some control points
that define the morphing of the domain. All the aforementioned
techniques are implemented in an open source Python package called PyGeM~\cite{pygem}. 
Earlier approaches to parameter space reduction
relies on modal analysis~\cite{forti2014efficient}, screening procedures based on Morris' randomized
one-at-a-time design~\cite{Morris,BallarinManzoniRozzaSalsa2013},
or semi-automatic reduction of the number of control points~\cite{DAmario}. The reduction of the
parameter space is achieved by retaining an optimal subset of the
possible control points. This can result in a set of admissible
deformations too shrunk. This can be overcome in part with the use of
active subspaces as explained in the following section. Another issue
is the optimal position of such control points. We underline that the
FFD is very versatile since it can be easily integrated into existing
software pipelines for CFD simulations, simply constructing a lattice of
points around the part of the domain to be morphed. In particular
PyGeM software can deal with a variety of file formats for both academic and
industrial purposes. 

\section{New Advances in Parameter Space Reduction with Active Subspaces}
\label{sec:Active_subspaces}
The improved capabilities in terms of computational power of the last
decade has led to more and more sophisticated CFD
simulations. Increasing the number of parameters allows finer
sensitivity analysis, and expand the design space to explore for shape
optimization problems. To fight the curse of dimensionality a possible
approach is to reduce the dimension of the parameter space. In the
last years the active subspaces technique, introduced
in~\cite{constantine2015active}, has been employed with success in
many different engineering problems such as optimal shape
design~\cite{lukaczyk2014active}, hydrologic
models~\cite{jefferson2015active}, naval
engineering~\cite{tezzele2017dimension,tezzele2018model}, and
uncertainty quantification~\cite{constantine2015exploiting}. Active
subspaces are a property of a parametric multivariate scalar
function, representing the quantity of interest, and a probability
density function used to sample the parameter space. An active
subspace is the span of the eigenvectors of the uncentered covariance
matrix of the gradients of the target function with respect to the
parameters. In practice the new reduced parameters are a linear
combination of the original ones, that accounts for how much the
quantity of interest varies along each parameter direction. It can
also be thought as a rotation of the parameter space in order to
unveil a lower dimensional behavior of the function of interest. The
insights given from the active subspace of a function are multiple: it
is possible to identify the more important parameters, and the ones we
can discard without affecting too much the approximation. This allows
an efficient choice of the geometrical parameters in optimal shape
design problems where the parametrization of the design is crucial
together with the exact choice of the parameters to describe it.

New advances have been made for what concerns time dependent
functions~\cite{constantine2017time}, the combination with POD-Galerkin
model reduction for cardiovascular
problems~\cite{tezzele2017combined}, the extension to vector-valued
functions~\cite{zahm2018gradient}, and the combination of different
active subspaces of different functions~\cite{ji2018shared}. Since,
to find an active subspace, we need couples of input/output data, the
technique can be easily integrated within existing pipelines
consisting, for example, on geometrical parametrization and
equation-free reduction techniques as in~\cite{tezzele2018model}.
Other approaches are possible to perform design-space nonlinear dimensionality
reduction, for a comparative review we suggest~\cite{van2009dimensionality}, while for an
application to naval engineering and shape optimization we cite~\cite{d2017nonlinear}.

\section{Naval and Nautical Engineering Applications}
\label{sec:Naval}

In the last decade, naval and nautical engineering fields have
witnessed a progressive introduction of high fidelity hydrodynamic
simulations into the design process of ship hulls, propellers and
other components. The increased computational capabilities and the
wealth of reliable simulation models and softwares, nowadays allows
for the evaluation of the fluid dynamic performance of virtual models
designed by the engineer. Such new virtual prototyping scenarios
propel the demand of new instruments which would make the design
pipeline more efficient. First, since many high fidelity CFD
simulations are typically carried out during a design campaign,
engineers now need reliable instruments to automatically produce a
suitable set of geometries which can be readily converted into quality
computational grids, and can at the same time explore in the most
extensive way the space of possible designs, in search for the optimal
one. The production of new shapes to be tested in fact, is in the most
common practice carried out by skilled designers who manually operate
CAD tools to obtain shapes that can improve hydrodynamic performance
and still fulfill the structural and bulk/volumetric constraints to
which the components are subjected. In addition to such time
consuming task, the geometries generated must be imported into mesh
generation tools to obtain suitable computational grids for CFD
simulations.  

Part of the current work at mathLab is focused on the development of
algorithms for the shape parametrization different components of the
ship considered in the design process. In such framework, both FFD
implemented in the PyGeM software and specifically developed shape
deformation strategies are applied to ship and planing yacht CAD
geometries, resulting in a series of parametrized IGES geometries. For
each morphed shape, we compute a series of geometrical and hydrostatic
parameters which are of typical interest of the naval and nautical
engineers. This allows for a first selection of the shapes that will
be tested, which avoid running detailed CFD simulation on
configurations that are \emph{a priori} known to be unsatisfactory. As
some mesh generators in the CFD community operate starting from STL
triangulations ---as do 3D printers--- a surface mesh generator has
been implemented to obtain water tight triangulations on non water
tight hull geometries. In addition, a further tool which employs PyGeM
FFD and RBF tools to directly deform the volumetric meshes generated
on the original hull shape is being developed. This will further
reduce the amount of human interaction required to carry out a
simulation campaign. As for the marine propeller shape
parametrization, a specific tool has been developed to build bottom-up
propeller virtual model in which the blades are parametrized on the
camber and thickness of the airfoils used, as well as the rake, skew,
pitch and chord distribution as a function of the radial
coordinate. The result of the procedure is represented by the IGES
geometry of the propeller, and the corresponding STL water-tight
triangulation. 

Along with the pre-processing tools described, the work is focusing
both on the development of reduced order models aimed at reducing the
computational cost of the hydrodynamics simulation campaign, and on a
smart post processing of the CFD simulations output based on active
subspaces. As for the first task, we are currently applying PODI tools
included in the ITHACA-FV to RANS ship hydrodynamic simulations
carried out with OpenFOAM. When a restricted amount of parameters are
considered, such model reduction strategy is in fact able to provide
reliable predictions of the entire flow field. For cases in which more
parameters have to be considered, more conventional POD approaches are
being considered. Once all the simulations have been carried out, the
relationship between input and output parameters is studied to
identify the presence of active subspaces. This information obtained
from such analysis is particularly relevant to design engineers. It in
fact is able to trigger possible redundancies in the parameter space
and, in some cases, it can also provide some physical insight on the
possible correlation among some parameters involved in the design
process. As mentioned, a first application of the entire pipeline
composed by shape parametrization, high fidelity computations and post
processing analysis with active subspaces has been presented
in~\cite{d2017nonlinear}.

\section{Industrial Engineering Applications}
\label{sec:Industrial}
\begin{figure}[t]
\centering
\includegraphics[width=0.7\textwidth]{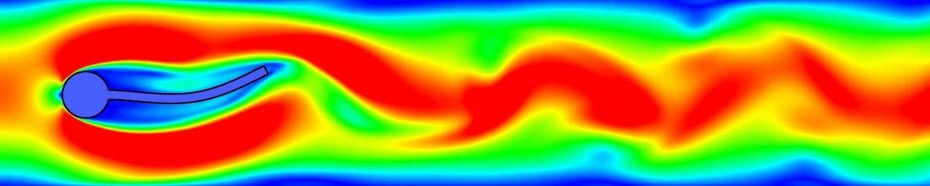}\\
\caption{Simulation of the Hron and Turek FSI benchmark \cite{Turek}.}
\label{fig:1}
\end{figure}

\begin{figure}[t]
\centering
\includegraphics[width=0.5\textwidth]{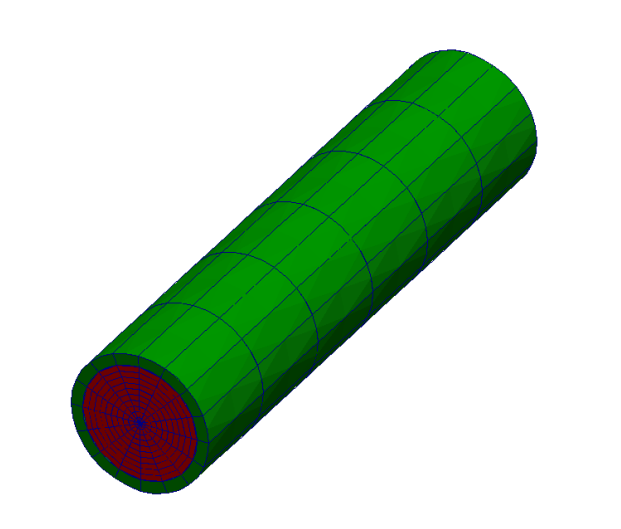}\\
\caption{Flexible tube benchmark problem with the employed fluid (red) and solid (green) computational grids \cite{Bogaers}.}
\label{fig:2}
\end{figure}

New simulation frameworks are required for the analysis of industrial
engineering problems that often involve multi-physics systems governed
by sets of coupled PDEs. A typical example belonging to this category
is the fluid structure interaction (FSI) that deals with the
investigation of the interaction phenomena between deformable or
movable structures with a surrounding (Figure \ref{fig:1}) or internal
(Figure \ref{fig:2}) fluid flow. The treatment of such problems
requires the introduction of proper coupling conditions as well as the
development of adequate numerical algorithms.  The scientific
literature includes a lot of works focused on the interaction between
an incompressible fluid flow and an hyperelastic solid but it should
be noted that industrial problems often are characterized by a very
complex multi-physics scenario which involves combined effects of
turbulent flows, thermo-chemical reactions, multiphase and interfacial
flows. Further elements of complexity that can be present in
industrial applications are for instance multiple regions separated by
multiple interfaces \cite{Bungartz}.

Regarding the full order model, the numerical procedures to solve FSI
problems may be classified into two categories: monolithic and
partitioned. In the monolithic approach, flow and structure equations
are solved simultaneously with a single solver and the interfacial
conditions are implicit in the solution procedure \cite{Hron}, whilst
in the partitioned approach flow and structure equations are  solved
in sequence  with  two  distinct solvers and the  interfacial
conditions  are  used  explicitly  to  communicate  information
between the fluid and structure solutions \cite{Schafer}. Of course
reduced order models should be considered related to both approaches. 

In the industrial context, the partitioned approach is the most
broadly used; it preserves software modularity because existing
solvers are coupled. Then, different and more efficient techniques
developed specifically for solving flow and structure equations can be
readily employed. Of course, the development of stable and accurate
coupling algorithms is required in partitioned approaches. In
particular, depending on the physics nature of the interaction,
one-way or two-way coupling procedures are demanded; in the former
case transfer quantities are sent from one domain to the other but not
in  the  opposite  direction,  in  the  latter  one  the  solver  data
is  always  transferred  both ways at the fluid-structure
interaction. Moreover, due to the strength of the coupling that can
occur in some problems as well as the well known problem of the
artificial added mass \cite{Forster}, that introduces further
instabilities, a partitioned approach often requires to iterate the
solution process of the systems of equations several times every time
step by determining a significant increase of the computational
cost. Quasi-Newton methods \cite{Degroote} are an example of efficient
iterative methods that have been employed in order to ensure a low
number of iterations. The numerical treatment of the interaction of a
fluid and slender structures is a very challenging problem
\cite{Paidoussis}.  Proper Orthogonal Decomposition (POD) reduced
order models have been built on a strong coupling algorithm for
partitioned FSI approaches in order to improve the convergence rate of
the iterative method when a lot of previous iterations are reused
\cite{Bogaers}. In \cite{Fang}, authors show that Singular Value
Decomposition (SVD) is also a very useful and robust algorithm to
treat the ill conditioning of the linear system involved by an
iterative method. A recent approach is presented in the study by
Ballarin et al \cite{Ballarin2017}.  

Monolithic approaches can potentially achieve better stability and
convergence properties but they require the development and handling
of a specialized code. Of course they are more robust for strongly
coupled problems but it has been shown that the monolithic approach
allows to obtain a good performance even for problems characterized by
a weak coupling \cite{Heil}. Some preliminary features of monolithic
reduced order models for FSI problems have been investigated in
\cite{Ballarin}.

\section{Conclusions}

With this manuscript, we have presented some examples of recent
advances in research activities of mathLab, SISSA group in various
domains focusing on the application of model order reduction
techniques. The recent developments in the area of model order
reduction methods are now on a level where there is better capability
to face much more complex problems. It is now possible to include data
driven modeling within the analysis, control and optimization. As
discussed in this manuscript, now model reduction techniques are
applied to several industrial engineering applications, as
demonstrated by the examples. Geometrical parameterization, shape
optimization and integration in the CAD modeling for the ships and
yachts design and analysis as well as geometrical reconstruction
through biomedical data demonstrate the strength of the methods
developed. As a last word, it looks even more promising, with the
growth in knowledge and experience in the field of model order
reduction, for computational scientists to be able to tackle more and
more challenging problems.

\section*{Acknowledgements}
The authors want to acknowledge the funding provided through the
European funding programme POR FESR 2014-2020, there are currently two
projects under POR FESR 2014-2020 program, namely Seakeeping of
Planning Hull of Yachts (SOPHYA) and Advanced methodologies for
hydro-acoustic design of naval propulsion (PRELICA). This work was
partially funded by the project HEaD, ``Higher
Education and Development'', supported by Regione FVG, European
Social Fund FSE 2014-2020, and by European Union
Funding for Research and Innovation --- Horizon 2020 Program --- in
the framework of European Research Council Executive Agency: H2020 ERC
CoG 2015 AROMA-CFD project 681447 ``Advanced Reduced Order Methods
with Applications in Computational Fluid Dynamics'' P.I. Gianluigi
Rozza.

\end{document}